\documentclass{amsart}

\usepackage{amssymb,latexsym,txfonts,amsmath,amsthm}
\usepackage{color,graphicx,pdfsync}
\usepackage[all,cmtip]{xy}

\usepackage{hyperref}
\hypersetup{
    colorlinks=true,       
    linkcolor=blue,          
    citecolor=blue,        
    filecolor=blue,      
    urlcolor=blue           
}

\newtheorem{theo}{Theorem}                     
\newtheorem{prop}{Proposition}                     
\newtheorem{coro}{Corollary}
\newtheorem{lemm}{Lemma}
\newtheorem{conj}{Conjecture}

\theoremstyle{remark}                  
\newtheorem{rema}{\bf Remark}

\title{A remark on the field of moduli of Riemann surfaces} 

\author{Rub\'en A. Hidalgo}

\address{
Departamento de Matem\'atica y Estad\'{\i}stica, Universidad de La Frontera. Temuco, Chile}
\email{ruben.hidalgo@ufrontera.cl}


\thanks{Supported by Project Fondecyt 1190001}

\begin{document}
\maketitle

\begin{abstract}
Let $S$ be a closed Riemann surface of genus $g\geq 2$ and let ${\rm Aut}(S)$ be its group of conformal automorphisms. It is well known that if either: (i) ${\rm Aut}(S)$ is trivial or (ii) $S/{\rm Aut}(S)$ is an orbifold of genus zero with exactly three cone points, then $S$ is definable over its field of moduli ${\mathcal M}(S)$. In the complementary situation, explicit examples for which ${\mathcal M}(S)$ is not a field of definition are known. We provide upper bounds for the minimal degree extension of ${\mathcal M}(S)$ by a field of definition in terms of the quotient orbifold $S/{\rm Aut}(S)$.
\end{abstract}

\section{Introduction}
There is a one-to-one correspondence between the categories of closed Riemann surfaces and non-singular irreducible complex projective algebraic curves (complex smooth curves for short). Two complex smooth curves $C_{1}$ and $C_{2}$ are isomorphic ($C_{1} \cong C_{2}$) if and only if the corresponding underlying Riemann surfaces structures are conformally equivalent. This equivalence permits to define the concepts of field of definition and field of moduli of a closed Riemann surface $S$ of genus $g$. A subfield ${\mathbb K}<{\mathbb C}$ is a field of definition of $S$
if it can be defined by a complex smooth curve given as the common zero locus of a  finite set of polynomials with coefficients in ${\mathbb K}$
(we also say that $S$ is definable over ${\mathbb K}$).  The field of moduli ${\mathcal M}(S)$ of $S$ is the intersection of all its fields of definitions. The conformal class of $S$, in the corresponding Riemann moduli space, is a ${\mathcal M}(S)$-rational point.  If $g=0$, then ${\mathcal M}(S)={\mathbb Q}$ (a model over ${\mathbb Q}$ is just the projective line ${\mathbb P}^{1}_{\mathbb Q}$). If $g=1$, then $S$ can be described by an elliptic curve of the form $E_{\lambda}:y^{2}z=x(x-z)(x-\lambda z)$, where $\lambda \in {\mathbb C}-\{0,1\}$. As $E_{\lambda_{1}}$ and $E_{\lambda_{2}}$ are isomorphic if and only if $j(\lambda_{1})=j(\lambda_{2})$, where 
$j$ is the classical Klein $j$-function, ${\mathcal M}(E_{\lambda})={\mathbb Q}(j(\lambda))$. It is known that $E_{\lambda}$ can be defined over ${\mathbb Q}(j(\lambda))$ \cite[Chapter III, Prop. 1.4]{Sil}. 

Let assume, from now on, that $g \geq 2$.
Explicit examples, not definable over their field of moduli, were provided by Earle \cite{Earle, Earle2}, Shimura \cite{Shimura} and Huggins \cite{Huggins} (in the hyperelliptic case) and by the author \cite{Hidalgo} and Kontogeorgis \cite{Kontogeorgis} (in the non-hyperelliptic situation). These examples are definable over a quadratic extension of their field of moduli. Necessary and sufficient conditions, for ${\mathcal M}(S)$ to be a field of definition, are given by 
Weil's Galois descent theorem \cite{Weil} (we recall it in Section \ref{Sec:prelim}).  
(If $S$ is definable over ${\mathcal M}(S)$, then the construction of a curve representing it, defined over ${\mathcal M}(S)$, is computationally hard. In \cite{HR} there is provided a computational method which permits to find such a model curve.) 
In \cite{Schwarz}, Schwarz  proved that the group  ${\rm Aut}(S)$ of conformal automorphisms of $S$ is finite. Later, in \cite{Hurwitz}, Hurwitz obtained the upper bound  $|{\rm Aut}(S)| \leq 84(g-1)$. 
If ${\rm Aut}(S)$ is trivial (the generic situation for $g \geq 3$), then Weil's conditions hold trivially, so $S$ can be defined over ${\mathcal M}(S)$. If the quotient orbifold $S/{\rm Aut}(S)$ has genus zero with exactly $3$ cone points (following Singerman, $S$ is called quasiplatonic), then Wolfart \cite{Wolfart} proved that $S$ can be defined over ${\mathcal M}(S)$. 
If $S$ is hyperelliptic, whose reduced group of automorphisms is neither trivial or cyclic, then Huggins \cite{Huggins-tesis,Huggins} proved that it is hyperelliptically definable over ${\mathcal M}(S)$. Recently, Lercier, Ritzenthaler and Sijsling \cite{LRS} have obtained that $S$ can be hyperelliptically defined over a quadratic extension of ${\mathcal M}(S)$ if its reduced group is a non-trivial cyclic group. In \cite{Mestre}, Mestre observed that if $g$ is even, then $S$ can be hyperelliptically defined over every field of definition of it (for $g$ odd this fact is in general false \cite{LR,LRS}). 

Let us denote by $({\rm FOD}/{\rm FOM})(S)$ the minimal degree extension of ${\mathcal M}(S)$ by a field of definition of $S$ (which is finite by \cite{Koizumi}).  As previously observed, for the known examples of Riemann surfaces $S$ not definable over their field of moduli,  $({\rm FOD}/{\rm FOM})(S)=2$. Our main result provides upper bounds for $({\rm FOD}/{\rm FOM})(S)$ in terms of the quotient orbifold $S/{\rm Aut}(S)$. 
Before to state such a result, we need to recall the definition of signature.
If $G<{\rm Aut}(S)$, then the quotient orbifold $S/G$ has an underlying Riemann surface structure $S_{G}$ of some genus $\gamma$ and there is a finite (possible empty) set of cone points, $p_{1},\ldots,p_{n}$. If $q_{j} \in S$ is projected to $p_{j}$, then the $G$-stabilizer of $q_{j}$ is a non-trivial cyclic group whose order $k_{j}$ is the cone order of $p_{j}$. We codify this information in a tuple $(\gamma;k_{1},\ldots,k_{n})$, called the signature of $S/G$. If $\gamma=0$, then we say that it is a triangular (respectively, quadrangular) signature if $n=3$ (respectively, $n=4$).

\begin{theo}\label{ejemplo1}
Let $S$ be a closed Riemann surface of genus $g \geq 2$.
{\rm (I)}  Assume that ${\rm Aut}(S)$ is non-trivial and let ${\mathcal O}=S/{\rm Aut}(S)$.
{\rm (1)}  If ${\mathcal O}$ is either of genus zero or hyperelliptic of even genus, then $({\rm FOD}/{\rm FOM})(S) \leq 2$.
{\rm (2)} If ${\mathcal O}$ is hyperelliptic of odd genus, then $({\rm FOD}/{\rm FOM})(S) \leq 4$.
{\rm (3)} If ${\mathcal O}$ has a signature of the form $(1;k_{1},\stackrel{r\geq 1}{\ldots} ,k_{1},k_{2},\ldots,k_{n})$, where $k_{1} \notin \{k_{2},\ldots,k_{n}\}$ ($k_{2},\ldots,k_{n}$ may have identical values), then (i) $({\rm FOD}/{\rm FOM})(S) \leq 2$,  if $r=1$, and (ii) $({\rm FOD}/{\rm FOM})(S) \leq r$, if $r \geq 2$.
{\rm (4)} If ${\mathcal O}$ is a non-hyperelliptic Riemann surface of genus $\gamma \geq 3$, then $({\rm FOD}/{\rm FOM})(S)\leq 2\gamma-2$.
{\rm (II)} $({\rm FOD}/{\rm FOM})(S)\leq 2(g-1)$.
\end{theo}

\begin{coro}\label{coro2}
If $S$ is a $n$-cyclic gonal closed Riemann surface ($n=2$ corresponds to hyperelliptic case), then $({\rm FOD}/{\rm FOM})(S) \leq 2$.
\end{coro}

Case (1) in Theorem \ref{ejemplo1} is well known (it was brought to the author's attention by C. Ritzenthaler). The main ingredient of the proof (which also holds in any algebraically closed field) is D\`ebes-Emsalem's theorem \cite{DE}. Part (II) is just a direct consequence of part (I) and the Riemann-Hurwitz formula. As an example, if $g \in \{2,3,4,5\}$, then Theorem \ref{ejemplo1} asserts that $S$ is definable over at most a quadratic extension of ${\mathcal M}(S)$, with possible exception at the following cases:
(1) $g=3$, ${\rm Aut}(S)=\langle \tau \rangle \cong {\mathbb Z}_{2}$ and $\#{\rm Fix}(\tau)=4$ ($({\rm FOD}/{\rm FOM})(S) \leq 4$); 
(2) $g=4$, ${\rm Aut}(S)=\langle \tau \rangle \cong {\mathbb Z}_{2}$ and $\#{\rm Fix}(\tau)=6$ ($({\rm FOD}/{\rm FOM})(S) \leq 6$); 
(3) $g=5$, ${\rm Aut}(S)=\langle \tau \rangle \cong {\mathbb Z}_{2}$ and $\#{\rm Fix}(\tau)\in\{0,8\}$ ($({\rm FOD}/{\rm FOM})(S) \leq 4$ in the first case and $({\rm FOD}/{\rm FOM})(S) \leq 8$ in the last one).
In Section \ref{Sec:quadrangular} we work out some particular examples of Riemann surfaces, called quadrangular.

\section{preliminaries}\label{Sec:prelim}
Let $C$ be a complex smooth curve given as the common zero locus of a finite set of complex polynomials $P_{1},\ldots,P_{r}$. If $\sigma \in {\rm Aut}({\mathbb C}/{\mathbb Q})$ (the group of field automorphisms of ${\mathbb C})$, then we let 
$C^{\sigma}$ be the algebraic curve defined as the common zero locus of the set of polynomials $P_{1}^{\sigma}, \ldots, P_{r}^{\sigma}$ (where $P_{j}^{\sigma}$ is the polynomial obtained from $P_{j}$ after applying $\sigma$ to its coefficients). The collection of 
those $\sigma$ such that  $C^{\sigma} \cong C$ form a subgroup $G_{C}<{\rm Aut}({\mathbb C}/{\mathbb Q})$ and its fixed field ${\mathcal M}(C)<{\mathbb C}$  is the field of moduli of $C$. It can be seen that $G_{C}={\rm Aut}({\mathbb C}/{\mathcal M}(C))$ and, by results due to Koizumi \cite{Koizumi}, that: (i) ${\mathcal M}(C)$ coincides with the intersection of all the fields of definitions of $C$ and (ii) there is a field of definition of $C$ that is a finite extension of ${\mathcal M}(C)$ (see also, \cite{HH}).

\subsection{Weil's Galois descent theorem}
Weil's Galois descent theorem provides necessary and sufficient conditions for an algebraic curve defined over any finite Galois extension ${\mathbb L}$ of a (perfect) field ${\mathbb K}$ to be definable also over ${\mathbb K}$.

\begin{theo}[Weil's Galois descent theorem \cite{Weil}]\
Let $C$ be a smooth curve, defined over a finite Galois extension ${\mathbb L}$ of a field ${\mathbb K}$.
If for every $\sigma \in {\rm Aut}({\mathbb L}/{\mathbb K})$ there is an isomorphism $f_{\sigma}:C \to C^{\sigma}$, defined over ${\mathbb L}$, such that, for all $\sigma, \tau \in {\rm Aut}({\mathbb L}/{\mathbb K})$ the compatibility condition $f_{\tau\sigma}=f^{\tau}_{\sigma} \circ f_{\tau}$
holds, then there exists a non-singular projective algebraic curve $E$, defined over ${\mathbb K}$, and there exists an isomorphism $R:C \to E$, defined over ${\mathbb L}$, such that $R^{\sigma} \circ f_{\sigma}=R$.
\end{theo}

In this paper, we are considering curves defined over ${\mathbb C}$; so a non-Galois extension of ${\mathbb Q}$. Next we observe how to deal in this case.
Let ${\mathbb K}$ be a subfield of ${\mathbb C}$ and $\overline{\mathbb K}$ be its algebraic closure in ${\mathbb C}$. Let $C$ be a complex smooth curve of genus $g\geq 2$ which is defined over $\overline{\mathbb K}$. As ${\rm Aut}(C)$ is finite, each of its elements is defined over $\overline{\mathbb K}$. Each $ \sigma \in {\rm Aut}({\mathbb C}/{\mathbb K})$ restricts to an element $\theta(\sigma) \in {\rm Aut}(\overline{\mathbb K}/{\mathbb K})$ and this provides a surjective homomorphism $\theta:{\rm Aut}({\mathbb C}/{\mathbb K}) \to {\rm Aut}(\overline{\mathbb K}/{\mathbb K})$ whose kernel is ${\rm Aut}({\mathbb C}/\overline{\mathbb K})$. Let us assume that for every $\sigma \in {\rm Aut}({\mathbb C}/{\mathbb K})$ there is an isomorphism $f_{\sigma}:C \to C^{\sigma}$. If $\tau \in {\rm Aut}({\mathbb C}/\overline{\mathbb K})$, then we get an isomorphism $f_{\sigma}^{\tau}:C=C^{\tau} \to C^{\sigma}=(C^{\sigma})^{\tau}$; so there is some $h_{\sigma} \in {\rm Aut}(C)$ with $f_{\sigma}^{\tau}=f_{\sigma} \circ h_{\sigma}$. Now, as ${\rm Aut}(C)$ is finite, it follows that the subgroup of ${\rm Aut}({\mathbb C}/\overline{\mathbb K})$ consisting of those elements $\tau$ with $f_{\sigma}^{\tau}=f_{\sigma}$ has finite index. In particular, $f_{\sigma}$ is defined over $\overline{\mathbb K}$. All the above asserts that it is possible to find a finite Galois extension ${\mathbb L}$ of ${\mathbb K}$ such that (i) $C$ and all its automorphisms are defined over ${\mathbb L}$ and (ii) for every $\sigma \in {\rm Aut}({\mathbb L}/{\mathbb K})$ there is an isomorphism $f_{\sigma}:C \to C^{\sigma}$ defined over ${\mathbb L}$. In this way, 
Weil's Galois descent theorem may be written, in our situation, as follows.

\begin{coro}\label{Prop:Weil}
Let $C$ be a complex smooth curve of genus $g \geq 2$ and let ${\mathbb K}$ be a subfield of ${\mathbb C}$ such that $C$ is defined over $\overline{\mathbb K}$.
If for every $\sigma \in {\rm Aut}({\mathbb C}/{\mathbb K})$ there is an isomorphism $f_{\sigma}:C \to C^{\sigma}$ such that, for all $\sigma, \tau \in {\rm Aut}({\mathbb C}/{\mathbb K})$, the compatibility condition $f_{\tau\sigma}=f^{\tau}_{\sigma} \circ f_{\tau}$
holds, then there exists a smooth curve $E$, defined over ${\mathbb K}$, and there exists an isomorphism $R:C \to E$, defined over $\overline{\mathbb K}$, such that $R^{\sigma} \circ f_{\sigma}=R$.
\end{coro}

As $C$ can be assume to be defined over a finite extension of ${\mathbb K}={\mathcal M}(C)$ \cite{Koizumi}, the above provides necessary and sufficient conditions for $C$ to be definable over ${\mathcal M}(C)$.

\subsection{D\`ebes-Emsalem's theorem}\label{Debes}
Next, we recall a consequence of Weil's Galois descent theorem (due to D\`ebes-Emsalem \cite{DE}) which will be useful to us. First, we need to state some notations and definitions. 
Let us consider complex smooth curves, say $C$ and $D$, and let $f:C \to D$ be a rational map which provides a (branched) holomorphic cover at the level of the corresponding closed Riemann surfaces. Let ${\mathbb K}$ be the field of moduli of $C$ and assume that $D$ is defined over it. By Koizumi's results \cite{Koizumi}, we may assume  that $C$ is defined over a finite Galois extension ${\mathbb L}$ of ${\mathbb K}$. 
We also assume that $C$ has genus at least two, so we may also assume that $f$ is defined over ${\mathbb L}$. Now, for each $\sigma \in {\rm Aut}({\mathbb C}/{\mathbb K})$ we consider the (branched) holomorphic cover $f^{\sigma}:C^{\sigma} \to D^{\sigma}=D$. 
We say that $f^{\sigma}:C^{\sigma} \to D$ and $f:C \to D$
are equivalent, noted as $\{f^{\sigma}:C^{\sigma} \to D\} \cong \{f:C \to D\}$, if there is an isomorphism $\phi_{\sigma}:C \to C^{\sigma}$ so that $f^{\sigma} \circ \phi_{\sigma}=f$. The field of moduli of the holomorphic cover $f:C \to D$, denoted by ${\mathcal M}(f:C \to D)$,  is the fixed field of the subgroup 
$\left\{\sigma \in {\rm Aut}({\mathbb C}/{\mathbb K}): \{f^{\sigma}:C^{\sigma} \to D\} \cong \{f:C \to D\}\right\}<{\rm Aut}({\mathbb C}/{\mathbb K}).$

\begin{theo}[D\`ebes-Emsalem \cite{DE}]\label{canonico}
Let $C$ be a complex smooth curve of genus $g \geq 2$ whose field of moduli is ${\mathbb K}$.
Then there exists a smooth curve $B$, defined over ${\mathbb K}$, and there exists 
a Galois holomorphic branched cover $Q:C \to B$, defined over $\overline{\mathbb K}$, with ${\rm Aut}(C)$ as its deck group, so that ${\mathcal M}(Q:C \to B)={\mathbb K}$. The curve $B$ is called the {\it ${\mathbb K}$-canonical model} of $C/{\rm Aut}(C)$.
Moreover, if $B$ contains at least one ${\mathbb L}$-rational point outside the branch locus of $Q$, where ${\mathbb L}$ is a finite extension of ${\mathbb K}$, then ${\mathbb L}$ is a field of definition of $C$. 
\end{theo}

\begin{rema}
Observe that the Galois branched covering $Q:C \to B$ can be defined over a given field if and only if $C$ itself can be.
\end{rema}

\section{Proof of Theorem \ref{ejemplo1}}
(I) Let $C$ be  a complex smooth curve, defined over a finite extension of its field of moduli ${\mathcal M}(S)={\mathbb K}$, defining $S$.
 Let $\overline{\mathbb K}$ be the algebraic closure of ${\mathbb K}$ in ${\mathbb C}$ and $\Gamma={\rm Aut}(\overline{\mathbb K}/{\mathbb K})$. We are assuming that ${\rm Aut}(C)$ is non-trivial. Following D\`ebes-Emsalem theorem, 
 there is an (smooth) algebraic curve $B$ (defined over ${\mathbb K}$) and a Galois holomorphic covering (defined over $\overline{\mathbb K}$), say $Q:C \to B$, whose deck group is ${\rm Aut}(C)$. If we are able to ensure the existence of a point on $B$, which is rational over an extension of ${\mathbb K}$ (outside the branch locus of $Q$), then (as a consequence of D\`ebes-Emsalem's theorem) $C$ will be definable over such an extension. 
(1) Assume $B$ has genus zero. The following idea was mentioned to the author by C. Ritzenthaler.
Let $K$ be the canonical divisor of $B$, defined over ${\mathbb K}$, and let us consider the Riemann-Roch's space 
${\mathcal L}(K)=\left\{f:B \to {\mathbb P}^{1}_{\overline{\mathbb K}}: (f)-K \geq 0\right\} \cup \{0\}$, which has dimension $3$ by Riemann-Roch's theorem. We may find a basis of ${\mathcal L}(K)$, say $f_{1}$, $f_{2}$ and $f_{3}$, where each element is defined over ${\mathbb K}$ (see Proposition 5.8. in Chapter 2 of \cite{Sil}). Construct, with such a basis, an isomorphism (defined over ${\mathbb K}$), say $L:B \to D \subset {\mathbb P}^{2}$, $L=[f_{1}:f_{2}:f_{3}]$, where $D$ is a (smooth conic or line) defined over ${\mathbb K}$. Now, on such a curve $D$ there are infinitely many rational points either over ${\mathbb K}$ or over a degree two extension of it. As $L$ is defined over ${\mathbb K}$, the pre-images of these points on $B$ are of the same type and we are done.
(2) Assume $B$ is hyperelliptic of even genus.
 Lemmas 4.2.1. and 4.2.2. in \cite{Huggins-tesis} assert the existence of an isomorphism $L:B \to W$ (defined over ${\mathbb K}$) where $W$ is the hyperelliptic curve defined by an equation of the form $y^{2}=F(x)$, where $F \in {\mathbb K}[x]$. As the curve $W$ has rational points over a suitable quadratic extension of ${\mathbb K}$, we may proceed similarly as above for the genus zero case. 
(3) Assume $B$ is hyperelliptic of odd genus $\gamma \geq 3$.
We follow similar arguments as in Section 2 in \cite{Mestre}. 
If $B$ were already hyperelliptically defined over ${\mathbb K}$, then it will have infinitely many rational points over a degree $2$ extension of ${\mathbb K}$ and we will be done. So, let us assume the model over ${\mathbb K}$ is not given in the hyperelliptic form. Take a base $\omega_{1}$,..., $\omega_{\gamma}$ of the space of holomorphic $1$-forms of $B$, each $\omega_{j}$ being ${\mathbb K}$-rational. Let us consider the two-to-one branched cover 
$P:B \to P(B)\subset {\mathbb P}_{\mathbb C}^{\gamma-1}:p \mapsto [\omega_{1}(p):\cdots:\omega_{\gamma}(p)].$
We know that $P(B)$ is a genus zero curve defined over ${\mathbb K}$ and that the image of a canonical divisor, defined over ${\mathbb K}$, of $B$ provides a ${\mathbb K}$-rational divisor of the form $2D$, where $D$ has degree $g-1$. The divisor $D$ is then ${\mathbb K}$-rational of even degree. This ensures the existence of infinitely many a rational points in $P(B)$ over an extension of degree $2$ of ${\mathbb K}$. By lifting such points to $B$ under $P$, we obtain infinitely many rational points over an extension of degree at most $4$ of ${\mathbb K}$ as desired.
(4) Assume $B$ has signature of the form $(1;k_{1},\stackrel{r \geq 1}{\ldots},k_{1},k_{2},\ldots,k_{n})$, where $k_{1} \notin \{k_{2},\ldots,k_{n}\}$. In this case, 
$B$ has genus one and there exactly $r\geq 1$ branch values of $Q$ of order $k_{1}$, say $p_{1},\ldots, p_{r} \in B$. If $r=1$, then we set $D=2p_{1} \in {\rm Div}(B)$. If $r \geq 2$, then we set  $D=p_{1}+\cdots+p_{r} \in {\rm Div}(B)$.
Note that, for every $\sigma \in \Gamma$, it holds that 
$D^{\sigma}=D$, that is, $D$ is ${\mathbb K}$-rational, and (i) ${\rm deg}(D)=2$, for $r=1$, and (ii) ${\rm deg}(D)=r$, for $r \geq 2$. By Riemann-Roch's theorem, ${\mathcal L}(-D)$ has dimension equal to ${\rm deg}(D)$. As $D$ is ${\mathbb K}$-rational, there is a basis for ${\mathcal L}(-D)$ formed by functions defined over ${\mathbb K}$; one of them is the constant $1$  and the others must be of order at most ${\rm deg}(D)$. Let $f$ be any of them and 
consider the map $f:B \to {\mathbb P}_{\mathbb K}^{1}$. As ${\mathbb P}_{\mathbb K}^{1}$ has infinitely many ${\mathbb K}$-rational points, $f$ is defined over ${\mathbb K}$ and it has degree at most ${\rm deg}(D)$, on $B$ there are rational points either over some extension of ${\mathbb K}$ of desired degree extension. 
(5) Assume $B$ is non-hyperelliptic of genus $\gamma \geq 3$.
Let us fix ${\mathbb K}$-rational holomorphic $1$-form, say $\omega_{1}, \ldots, \omega_{\gamma}$, defining a basis for the space of holomorphic $1$-forms of $B$, and  consider the holomorphic canonical embedding (which is ${\mathbb K}$-rational) $\phi:B \to {\mathbb P}^{\gamma-1}: p \mapsto [\omega_{1}(p):\cdots:\omega_{\gamma}(p)]$. As ${\mathcal Q}=\phi(B)$ is a smooth curve of degree $2\gamma-2$, defined over ${\mathbb K}$, it has infinitely many rational points over an extension of degree at most $2\gamma-2$ of ${\mathbb K}$. As $\phi$ has degree one and is defined over ${\mathbb K}$, the same holds for $B$.
(II) As $S$ is definable over ${\mathcal M}(S)$ if ${\rm Aut}(S)$ is trivial, we assume this group to be non-trivial.
The case $g=2$ follows from part (1). If $g \geq 3$, as $2(g-1)\geq 4$, then we only need to assume that $S/{\rm Aut}(S)$ is either (a) non-hyperelliptic of genus $\gamma \geq 3$ or (b) of signature $(1;k_{1},\stackrel{r \geq 2}{\ldots},k_{1},k_{2},\ldots,k_{n})$, where $k_{1} \notin \{k_{2},\ldots,k_{n}\}$. In case (a), as $g \geq \gamma$, it follows that $g-1 \geq \gamma-1$ and the result follows from part (4). In case (b), it follows from 
the Riemann-Hurwitz formula that $2(g-1) \geq r|{\rm Aut}(S)|(1-1/k_{1})\geq r$ and the result follows from part (3).

\section{Example: Quadrangular surfaces}\label{Sec:quadrangular}
A Riemann surface $S$ of genus $g \geq 2$ is called quadrangular if there exists $H<{\rm Aut}(S)$ with $S/H$ of quadrangular signature. In this case, by Theorem \ref{ejemplo1}, $({\rm FOD}/{\rm FOM})(S)\leq 2$. 

\begin{lemm}\label{lema1}
If $S$ is quadrangular, then the signature of 
$S/{\rm Aut}(S)$ is either triangular or quadrangular.
\end{lemm}
\begin{proof}
This is clear if $H={\rm Aut}(S)$. Let us assume that $H$ has index $N \geq 2$ in ${\rm Aut}(S)$. Let $(0;k_{1},k_{2},k_{3},k_{4})$ be the signature of $S/H$ and $(\gamma;p_{1},\ldots,p_{m})$, $p_{j} \geq 2$ for $m>0$, be the signature of $S/{\rm Aut}(S)$. The hyperbolic area of these hyperbolic orbifolds are respectively (see, for instance, \cite{Beardon}):
${\rm Area}(S/{\rm Aut}(S))=2 \pi \left(2\gamma-2 + \sum_{j=1}^{m}(1-p_{j}^{-1})\right)$ and 
${\rm Area}(S/H)=2 \pi \left( -2 + \sum_{j=1}^{3}(1-k_{j}^{-1})\right).$
By the Riemann-Hurwitz formula,  $|{\rm Aut}(S)| {\rm Area}(S/{\rm Aut}(S))={\rm Area}(S)=|H| {\rm Area}(S/H)$, so $2{\rm Area}(S/{\rm Aut}(S)) \leq N{\rm Area}(S/{\rm Aut}(S)) \leq {\rm Area}(S/H)$, that is,
$$(*) \quad 0<2\gamma-2+\sum_{j=1}^{m}\left(1-p_{j}^{-1}\right) \leq  \left(1- \sum_{j=1}^{3} k_{j}^{-1}\right) /2<1/2.$$

In particular, $2\gamma-2<1/2$, that is, $\gamma \in \{0,1\}$. If $\gamma=1$, then $\sum_{j=1}^{m}(1-p_{j}^{-1})<1/2$, a contradiction (as $p_{j} \geq 2$). So, $\gamma=0$ (and $m \geq 3$). In this way, $(*)$ implies $-2+\sum_{j=1}^{m}(1-p_{j}^{-1})<1/2$, and (as $p_{j} \geq 2$), this asserts that $-2+m/2<1/2$, that is, $m \leq 4$.
\end{proof}

\subsection{Quasiplatonic quadrangular surfaces}
Examples of quadrangular surfaces, called quasiplatonic quadrangular, are those for which $S/{\rm Aut}(S)$ has a quadrangular signature with one of the cone orders different from the orders  three (they are definable over their field of moduli \cite{AQ}). 

\begin{prop}\label{coro1}
Let $S$ be such that 
$S/{\rm Aut}(S)$ has signature $(0;a,b,c,d)$, where $d \notin \{a,b,c\}$. Assume that the cone points of $S/{\rm Aut}(S)$ (which we may identify with $\widehat{\mathbb C}$) are $\infty$, $0$, $1$ and $\lambda$, where $\lambda$ has the order $d$. Then ${\mathcal M}(S)$ is an extension of degree $e \in \{1,2,3,4,6,9,12,18, 36\}$ of ${\mathbb Q}(j(\lambda))$, where $j$ is the Klein $j$-function
$j(\lambda)=(1-\lambda+\lambda^{2})^{3}/\lambda^{2}(1-\lambda)^{2}$. 
\end{prop}
\begin{proof}
Let $P:S \to \widehat{\mathbb C}$ be a regular branched cover with ${\rm Aut}(S)$ as its deck group whose branch values are $\infty$, $0$, $1$ and $\lambda \in {\mathbb C}-\{0,1\}$, such that $a$ is the order of $\infty$, $b$ is the order of $0$, $c$ is the order of $1$, $d$ is the order of $\lambda$. If $\sigma \in G_{S}$, then there is an isomorphism $f_{\sigma}:S \to S^{\sigma}$ and a M\"obius transformation $M_{\sigma}$ such that $M_{\sigma} \circ P=P^{\sigma} \circ f_{\sigma}$. As $M_{\sigma}$  sends the set of branch values of $P$ onto the set of branch values of $P^{\sigma}$ and preserves the branch order, it follows that (i) $M_{\sigma}(\lambda)=\sigma(\lambda)$ and (ii) $M_{\sigma}$ preserves the set $\{\infty,0,1\}$. By (ii), $M_{\sigma} \in J =\langle A(z)=1/z, B(z)=z/z-1 \rangle \cong {\mathfrak S}_{3}$. Now, by (i), $\sigma(\lambda)$ belongs to the $J$-orbit of $\lambda$, that is, $j(\lambda)=j(\sigma(\lambda))=\sigma(j(\lambda))$, in particular, ${\mathbb Q}(j(\lambda))<{\mathcal M}(S)$. The collection $\{M_{\sigma}: \sigma \in G_{S}\}$ satisfies the co-cycle conditions in Weil's Galois descent theorem \cite{Weil}, that is, $M_{\tau\sigma}=M_{\sigma}^{\tau} \circ M_{\tau}$, for every $\sigma, \tau \in G_{S}$. This fact, together with $M_{\sigma}^{\tau}=M_{\sigma}$ (by (ii)), asserts that $\Theta:G_{S} \to J: \sigma \mapsto M_{\sigma}^{-1}$, is an homomorphism. As the kernel of $\Theta$ consists of those $\sigma \in G_{S}$ such that $\sigma(\lambda)=\lambda$, it follows that ${\mathcal M}(S)$ is an extension of ${\mathcal M}(S) \cap {\mathbb Q}(\lambda)$ of degree $d \in \{1,2,3,6\}$. As ${\mathbb Q}(j(\lambda))<{\mathcal M}(S) \cap {\mathbb Q}(\lambda)<{\mathbb Q}(\lambda)$, it also follows that ${\mathcal M}(S) \cap {\mathbb Q}(\lambda)$ is an extension of degree $s \in \{1,2,3,6\}$ of ${\mathbb Q}(j(\lambda))$.

\end{proof}

\subsubsection{\bf Example: genus three quasiplatonic quadrangular surfaces}
In \cite{Fuertes-Streit}, Fuertes and Streit described the quadrangular genus three surfaces $S$. In this case, $S/{\rm Aut}(S)$ has as signature one of the following: $(0;2,2,2,3)$, $(0;2,2,2,4)$, $(0;2,2,2,6)$ and $(0;2,3,3,6)$. In particular, they are quasiplatonic quadrangular, so $({\rm FOD}/{\rm FOM})(S)=1$. The automorphism groups of these surfaces were also discussed in \cite{LRRS,LRS}. These surfaces are describe below.

\noindent
{\bf (I)} ${\mathcal F}_{1}: \quad C_{\lambda}=\{x^{4}+y^{4}+z^{4}+\lambda (x^{2}y^{2}+y^{2}z^{2}+z^{2}x^{2})=0\}$, $\lambda \in {\mathcal P}_{1}:={\mathbb C}-\{\pm -1, \pm 2, -3(1\pm \sqrt{-7})/2 \}$. In this case, ${\rm Aut}(C_{\lambda}) \cong {\mathfrak S}_{4}$ and the signature is $(0;2,2,2,3)$. In  \cite{KFT} the family ${\mathcal F}_{1}$ is also called the KFT family. If $\lambda_{1}, \lambda_{2} \in {\mathcal P}_{1}$, then $C_{\lambda_{1}}$ and $C_{\lambda_{2}}$ are conformally equivalent if and only if  $\lambda_{1}=\lambda_{2}$. In particular, ${\mathcal M}(C_{\lambda})={\mathbb Q}(\lambda)$.

\noindent
{\bf (II)} ${\mathcal F}_{2}: \quad B_{\lambda}=\{y^{4}=x^{4}-\lambda x^{2}z^{2}+z^{4}\}$, $\lambda \in {\mathcal P}_{2} \subset {\mathbb C}-\{ \pm 2 \}$.
In this case,  ${\rm Aut}(B_{\lambda}) \cong {\mathbb Z}_{2} \ltimes ({\mathbb Z}_{2} \times {\mathbb Z}_{4})$ and the signature is $(0;2,2,2,4)$. If $\lambda_{1}, \lambda_{2} \in {\mathcal P}_{2}$, then $B_{\lambda_{1}}$ and $B_{\lambda_{2}}$ are conformally equivalent if and only if  $\lambda_{1}^{2}(\lambda_{1}^{2}-36)^{2}(\lambda_{2}^{2}-4)^{2}= \lambda_{2}^{2}(\lambda_{2}^{2}-36)^{2}(\lambda_{1}^{2}-4)^{2}$. In particular, ${\mathcal M}(B_{\lambda})={\mathbb Q}(\frac{\lambda^{2}(\lambda^{2}-36)^{2}}{(\lambda^{2}-4)^{2}})$.

\noindent
{\bf (III)} ${\mathcal F}_{3}: \quad D_{\lambda}=\{y^{2}z^{6}=x^{8}-\lambda x^{4}z^{4}+z^{8}\}$, $\lambda \in {\mathcal P}_{3} \subset {\mathbb C}-\{ \pm 2 \}$. In this case, ${\rm Aut}(D_{\lambda}) \cong {\mathbb Z}_{2} \times {\mathbb D}_{4}$ and the signature is $(0;2,2,2,4)$. If $\lambda_{1}, \lambda_{2} \in {\mathcal P}_{3}$, then $D_{\lambda_{1}}$ and $D_{\lambda_{2}}$ are conformally equivalent if and only if  $\lambda_{1}^{2}=\lambda_{2}^{2}$. In particular, ${\mathcal M}(D_{\lambda})={\mathbb Q}(\lambda^{2})$. If $\lambda \neq 0$, by making the change of variables $v=x$, $u=\sqrt[8]{\lambda^{6}} \; y$, and $w=\left(\sqrt[4]{\lambda}\right)^{-1}  z$, 
the curve $B_{\lambda}$ is isomorphic to the curve $\{u^{2}w^{6}=x^{8}+\lambda^{2}w^{4}(w^{4}-x^{4})\} \subset {\mathbb P}_{\mathbb C}^{2}$.

\noindent
{\bf (IV)} ${\mathcal F}_{4}: \quad E_{\lambda}=\{y^{2}z^{5}=x^{7}-\lambda x^{4}z^{3}+xz^{6}\}$, $\lambda \in {\mathcal P}_{4} \subset {\mathbb C}-\{ \pm 2\}$. In this case, ${\rm Aut}(E_{\lambda}) \cong {\mathbb D}_{6}$ and the signature is $(0;2,2,2,6)$. If $\lambda_{1}, \lambda_{2} \in {\mathcal P}_{4}$, then $E_{\lambda_{1}}$ and $E_{\lambda_{2}}$ are conformally equivalent if and only if  $\lambda_{1}^{2}=\lambda_{2}^{2}$. In particular, ${\mathcal M}(E_{\lambda})={\mathbb Q}(\lambda^{2})$. If $\lambda \neq 0$, by making the change of variables
$v=x$,  $u=\sqrt[6]{\lambda^{5}}\; y$, and  $w=\left(\sqrt[3]{\lambda}\right)^{-1}  z$,
the curve $E_{\lambda}$ is isomorphic to the curve 
$\{u^{2}w^{5}=x^{7}+\lambda^{2}x w^{3}(w^{3}-x^{3})\} \subset {\mathbb P}_{\mathbb C}^{2}$.

\noindent
{\bf (V)} ${\mathcal F}_{5}: \quad F_{\lambda}=\{y^{3}z=x^{4}-(1+\lambda)x^{2}z^{2}+\lambda z^{4}\}$, $\lambda \in {\mathcal P}_{5} \subset {\mathbb C}-\{ 0,1 \}$. In this case, ${\rm Aut}(F_{\lambda}) \cong {\mathbb Z}_{6}$ and the signature is $(0;2,3,3,6)$.  If $\lambda_{1}, \lambda_{2} \in {\mathcal P}_{5}$, then $F_{\lambda_{1}}$ and $F_{\lambda_{2}}$ are conformally equivalent if and only if  $\lambda_{1}+\lambda_{1}^{-1}=\lambda_{2}+\lambda_{2}^{-1}$. In particular, ${\mathcal M}(F_{\lambda})={\mathbb Q}(\lambda+\lambda^{-1})$. If $\lambda \neq -1$, by making the change of variables 
$v=\left(\frac{\sqrt{1+\lambda}}{\sqrt[4]{\lambda}}\right)x$,  
$u=\sqrt[12]{\lambda} \; y$,  and 
$w=\sqrt[4]{\lambda} \; z$,
 $F_{\lambda}$ is  isomorphic to 
$\left\{u^{3}w=\left(2+ \lambda+\lambda^{-1}  \right)^{-1} t^{4}-t^{2}w^{2}+w^{4}\right\} \subset {\mathbb P}_{\mathbb C}^{2}$.

\begin{rema}[Case of genus four]
In \cite{Swinarski}, Swinarski computed explicit algebraic curves for genus four Riemann surfaces $S$ such that $S/{\rm Aut}(S)$ has quadrangular signature. These possible signatures are: 
$(0;2,2,2,3)$, $(0;2,2,2,4)$, $(0;2,2,2,5)$, $(0;2,2,2,8)$ and $(0;2,2,3,3)$. With the exception of the last signature, $S$ is quasiplatonic quadrangular.
\end{rema}

\subsection{Quadrangular homology curves}
Let $S$ be a closed Riemann surface of genus $g \geq 2$ and $H<{\rm Aut}(S)$ be an abelian group. We say that $H$ is a homology group of $S$ if:
(i) $S/H$ has genus $0$ and (ii) there is no other closed Riemann surface $R$, of genus greater than of $S$, admitting an abelian group $K<{\rm Aut}(R)$ such that $R/K$ and $S/H$ are isomorphic Riemann orbifolds. Equivalently, if $\Gamma$ is a Fuchsian group such that the orbifolds  
${\mathbb H}^{2}/\Gamma$ and $S/H$ are isomorphic as Riemann orbifolds, then $H$ is a homology group if the derived subgroup $\Gamma'$ is torsion free and there is a biholomorphism $\phi:{\mathbb H}^{2}/\Gamma' \to S$ conjugating $\Gamma/\Gamma'$ to $H$. In this case, let 
${\rm Aut}_{H}(S)$ be the normalizer of $H$ in ${\rm Aut}(S)$ and let ${\rm Aut}(S/H)$ be the group of conformal automorphisms of the orbifold $S/H$ (this is exactly the subgroup of M\"obius transformations keeping invariant the set of cone points of $S/H$ and preserving their cone orders). If $\pi:S \to S/H$ is a regular branched covering with deck group $H$, then there is a natural homomorphism $\theta:{\rm Aut}_{H}(S) \to {\rm Aut}(S/H)$ such that $\pi \circ \alpha = \theta(\alpha) \circ \pi$, for $\alpha \in {\rm Aut}_{H}(S)$, whose kernel is $H$. Let $\beta \in {\rm Aut}(S/H)$ and let $\widehat{\beta} \in {\rm Aut}({\mathbb H}^{2})$ be a lifting of it. As $\widehat{\beta}$ normalizes $\Gamma$ and $\Gamma'$ is a characteristic subgroup of $\Gamma$, $\widehat{\beta}$ also normalizes $\Gamma'$ and, in particular, it induces a conformal automorphism $\alpha \in {\rm Aut}(S)$ such that $\theta(\alpha)=\beta$. So, there is a natural short exact sequence
\begin{equation}\label{exacto}
1 \to H \to {\rm Aut}_{H}(S) \stackrel{\theta}{\to} {\rm Aut}(S/H) \to 1.
\end{equation}

A closed Riemann surface $S$, of genus $g \geq 2$, is called a homology Riemann surface if it admits an homology group $H<{\rm Aut}(S)$. If moreover, $S/H$ has signature $(0;n_{1},...,n_{4})$, then we say that $S$ is a quadrangular homology Riemann surface.

\begin{rema}
(1) If $H$ is a homology group with $S/H$ of signature $(0;k,\stackrel{n+1}{\ldots},k)$, then $H \cong {\mathbb Z}_{k}^{n}$. If $(k-1)(n-1)>2$, then 
in \cite{HLKP} it was observed that $H$ is unique, so a normal subgroup (the case $n=3$ was obtained in \cite{FGHL}). 
(2) A homology Riemann surface may admit different homology groups. For instance, consider the hyperelliptic Riemann surface defined by the hyperelliptic curve $y^{2}=x^{2g+2}-1$, where $g \geq 2$ is even. Its hyperelliptic involution is $\tau(x,y)=(x,-y)$ and $S$ admits the order $2g+2$ automorphism
$\alpha(x,y)=(e^{\pi i/(g+1)}x,y)$. If $A=\langle \alpha, \tau \rangle \cong {\mathbb Z}_{2g+2} \times {\mathbb Z}_{2}$, then one may see that the quotient orbifold $S/A$ has signature $(0;2,2g+2,2g+2)$ and $A$ is a homology group of $S$. Similarly, if $B=\langle \alpha^{2}, \tau \rangle \cong {\mathbb Z}_{2g+2}$, then $S/B$ has signature $(0;2,2,g+1,g+1)$ and $B$ is again a homology group of $S$.  
\end{rema}

\begin{theo}\label{Main}
If $S$ is a quadrangular homology Riemann surface, then $({\rm FOD}/{\rm FOM})(S)=1$.
\end{theo}
\begin{proof}
Let $H<{\rm Aut}(S)$ be a homology group of $S$ with $S/H={\mathcal O}$ of signature $(0;a,b,c,d)$. By Lemma \ref{lema1}, either $S$ is quasiplatonic (so it can be defined over their fields of moduli  \cite{Wolfart}) or $S/{\rm Aut}(S)$ has signature of the form $(0;p_{1},p_{2},p_{3},p_{4})$. Let us assume we are in the last case. By the results in \cite{AQ}, we only need to take care of either (i) $p_{1}=p_{2}=p_{3}=p_{4}\geq 3$ or (iii) $p_{1}=p_{2} \neq p_{3}=p_{4}$ (up to permutation of indices). As in any of these two situations ${\rm Aut}(S/H)$ is non-trivial, it follows from the short exact sequence \eqref{exacto} that $H \neq {\rm Aut}(S)$. There are Fuchsian groups $\Gamma_{1}<\Gamma_{2}$, $\Gamma_{1}$ of fnite index in $\Gamma_{2}$, such that ${\mathbb H}^{2}/\Gamma_{1}=S/H$, ${\mathbb H}^{2}/\Gamma_{2}=S/{\rm Aut}(S)$. As both groups have quadrangular signature, the Teichm\"uller spaces $T(\Gamma_{1})$ and $T(\Gamma_{2})$ have the same real dimension $2$. It follows from 
(Theorems 1 and Proposition 5 in \cite{Singerman}) that $\Gamma_{2}$ is a normal subgroup of $\Gamma_{1}$ and (up to permutation of the elements in both quadruple signatures) either:
(1) $(a,b,c,d)=(k,k,k,k)$ and $(p_{1},p_{2},p_{3},p_{4})=(2,2,2,t)$; or 
(2) $(a,b,c,d)=(t_{1},t_{1},t_{2},t_{2})$ and $(p_{1},p_{2},p_{3},p_{4})=(2,2,t_{1},t_{2})$.
In case (1), as $g\geq 2$, it follows that $t \geq 3$ and the result follows from \cite{AQ} (similarly for case (2) if $t_{1} \neq t_{2}$). So, let us assume we are in case (2) with $t_{1}=t_{2}=k \geq 3$, that is, 
$(a,b,c,d)=(k,k,k,k)$ and $(p_{1},p_{2},p_{3},p_{4})=(2,2,k,k)$. We claim that this last case is not possible. By \cite{FGHL,HLKP}, 
$H$ is a normal subgroup of ${\rm Aut}(S)$, i.e., ${\rm Aut}_{H}(S)={\rm Aut}(S)$ and we have the short exact sequence
$1 \to H \to {\rm Aut}(S) \to {\rm Aut}({\mathcal O}) \to 1,$
that is, $S/{\rm Aut}(S)={\mathcal O}/{\rm Aut}({\mathcal O})$ and that ${\rm Aut}(S)/H={\rm Aut}({\mathcal O})$. Notice that ${\rm Aut}({\mathcal O})$ contains as subgroup the Klein group
$\widehat{J}=\left\{I(z)=z, A(z)=\lambda/z, B(z)=(z-\lambda)/(z-1), C(z)=\lambda(1-z)/(\lambda-z)\right\} \cong {\mathbb Z}_{2}^{2}.$
It follows the existence of a subgroup $K<{\rm Aut}(S)$ such that $H \lhd K$ and $K/H=\widehat{J}$. So $S/K=
{\mathcal O}/\widehat{J}$ has signature $(0;2,2,2,k)$, with $k \geq 3$. Now, following similar computations as done in the proof of Lemma \ref{lema1} (also by \cite{Singerman}), it can be obtained that the signature of $S/{\rm Aut}(S)$ cannot be of the form $(0;2,2,k,k)$.
\end{proof}

\begin{rema}
An homology Riemann surface of genus $g\geq 2$ and type $(k,k,k,k)$ (also called a generalized Fermat curve of type $(k,3)$), with $k \geq 3$, can be represented by a curve
$C_{\lambda}=\left\{ \begin{array}{c}
x_{1}^{k}+x_{2}^{k}+x_{3}^{k}=0=
\lambda x_{1}^{k}+x_{2}^{k}+x_{4}^{k}
\end{array}
\right\}\subset {\mathbb P}_{\mathbb C}^{3}$, where $\lambda \in {\mathbb C}-\{0,1\}$. By Theorem \ref{Main}, $C_{\lambda}$ is definable over its field of moduli.
In \cite{HJ} it was observed that ${\mathcal M}(C_{\lambda})={\mathbb Q}(j(\lambda))$ (note the similarity with the case of elliptic curves). 
\end{rema}

\begin{conj}
If there an abelian group $H<{\rm Aut}(S)$ such that $S/H$ has signature of the form $(0;a,b,c,d)$, then $S$ can be defined over its field of moduli.
\end{conj}

\subsection{Hyperelliptic quadrangular Riemann surfaces}
If $S$ is a hyperelliptic Riemann surface, then 
Theorem \ref{ejemplo1} asserts that $({\rm FOD}/{\rm FOM})(S)\leq 2$. If, moreover, $S$ is quadrangular, then it is definable over ${\mathcal M}(S)$.

\begin{coro}\label{hipereliptica}
If $S$ is hyperelliptic quadrangular, then $({\rm FOD}/{\rm FOM})(S)=1$.
\end{coro}
\begin{proof}
Let $\iota:S \to S$ be the the hyperelliptic involution of $S$ and let $H<{\rm Aut}(S)$ such that $S/H$ has quadrangular signature $(0;a,b,c,d)$. As observed at the beginning of Section \ref{Sec:quadrangular}, either (i) $S$ is quasiplatonic or (ii) the signature of $S/{\rm Aut}(S)$ is of type $(0;r,s,t,u)$.  
In case (i), by \cite{Wolfart}, $({\rm FOD}/{\rm FOM})(S)=1$. Let us now consider case (ii).
As $S/\langle \iota \rangle$ is an orbifold with signature $(0;2,\stackrel{2g+2}{\cdots},2)$ and $2g+2>4$, necessarily ${\rm Aut}(S) \neq \langle \iota \rangle$. If $g=2$, then, by \cite{Cardona},  $({\rm FOD}/{\rm FOM})(S)=1$. Let us assume now on that $g \geq 3$. Since $({\rm FOD}/{\rm FOM})(S)=1$ 
when ${\rm Aut}(S)/\langle \iota \rangle$ is not a cyclic group \cite{Huggins}, we only need to assume ${\rm Aut}(S)/\langle \iota \rangle \cong {\mathbb Z}_{n}$, for a suitable $n \geq 2$. Let $P:S \to \widehat{\mathbb C}$ be a regular two-fold branched cover with deck group $\langle \iota\rangle$. Up to composition of $P$ at the left by a suitable M\"obius transformation, we may assume that ${\rm Aut}(S)/\langle \iota \rangle=\langle A(z)=e^{2 \pi i/n}z \rangle \cong {\mathbb Z}_{n}$. As $S/{\rm Aut}(S)$ has a quadrangular signature and the set of branch values of $P$ is invariant under the action of $A$, it follows that the set of branch values of $P$ can only be of one of the following forms (up to composition of $P$ at the left by a M\"obius transformation of the form $T(z)=qz$ or $T(z)=q/z$, for a suitable $q \neq 0$, where $\omega=e^{2 \pi i/n}$):
(1) $\{\omega^{k}, \mu\omega^{k}: k=1,...,n\}$, where $\mu \in {\mathbb C}$ with $\mu^{n} \neq 1$; or
(2) $\{\infty\} \cup \{\omega^{k}, \mu\omega^{k}: k=1,...,n\}$, where $\mu \in {\mathbb C}$ with $\mu^{n} \neq 1$; or
(3) $\{\infty,0\} \cup \{\omega^{k}, \mu\omega^{k}: k=1,...,n\}$, where $\mu \in {\mathbb C}$ with $\mu^{n} \neq 1$.
As the number of branch values of $P$ is even, (2) is not possible. In cases (1) and (3) we have that the M\"obius transformation $B(z)=\mu/z$ keeps invariant the branch values of $P$. It follows that ${\rm Aut}(S)/\langle \iota \rangle$ contains the dihedral group $D_{n}=\langle A,B\rangle$, a contradiction.
\end{proof}

\subsection{Acknowledgements}
The author would like to express his deep gratitude to the referee for supplying very useful comments, suggestions and corrections. 



\begin{thebibliography}{99}

\bibitem{AQ}
M. Artebani and S. Quispe. 
Fields of moduli and fields of definition of odd signature curves. 
{\it Arch. Math.} (Basel) {\bf 99} No. 4 (2012), 333--344.


\bibitem{Beardon}
A. F. Beardon.
{\it The Geometry of Discrete Groups}.
Graduate Texts in Mathematics {\bf 91}, Springer-Verlag New York, 1983.



\bibitem{DE}
P. D\`ebes and M. Emsalem.
On Fields of Moduli of Curves.
{\it J. of Algebra} {\bf 211} (1999) 42--56.


\bibitem{Cardona}
G. Cardona and J. Quer.
Field of moduli and field of definition for curves of genus 2. Computational aspects of algebraic curves, 71--83. 
{\it Lecture Notes Ser. Comput.} {\bf 13}, World Sci. Publ., Hackensack, NJ, 2005.


\bibitem{Earle}
C. J. Earle.
On the moduli of closed Riemann surfaces with symmetries.
{\it Advances in the Theory of Riemann Surfaces} (1971) 119-130. Ed. L.V. Ahlfors et al. 
(Princeton Univ. Press, Princeton).


\bibitem{Earle2}
C. J. Earle.
Diffeomorphisms and automorphisms of compact hyperbolic 2-orbifolds. 
In {\it Geometry of Riemann surfaces}, {\it London Math. Soc. Lecture Note Ser.} {\bf 368}, 139--155.
Cambridge Univ. Press, Cambridge, 2010.



\bibitem{FGHL}
Y. Fuertes, G. Gonz\'alez-Diez, R. A. Hidalgo and M. Leyton.
Automorphism group of Generalized Fermat curves of type  $(k,3)$.
{\it Journal of Pure and Applied Algebra} {\bf 217} No. 10 (2013), 1791--1806.


\bibitem{Fuertes-Streit}
Y. Fuertes and M. Streit.
Genus $3$ normal coverings of the Riemann sphere branched over $4$ points.
{\it Rev. Mat. Iberoamericana} {\bf 22} No. 2 (2006), 413--454.


\bibitem{HH}
H. Hammer and F. Herrlich.
A Remark on the Moduli Field of a Curve.
{\it Arch. Math.} {\bf 81} (2003) 5--10.

\bibitem{Hidalgo}
R. A. Hidalgo.
Non-hyperelliptic Riemann surfaces with real field of moduli but not definable over the reals. 
{\it Archiv der Mathematik} {\bf 93} (2009) 219--222.

\bibitem{HLKP}
R. A. Hidalgo, A. Kontogeorgis, M. Leyton-{\'A}lvarez and P. Paramantzoglou.
Automorphisms of the Generalized Fermat curves.
{\it Journal of Pure and Applied Algebra} {\bf 221} (2017), 2312--2337.

\bibitem{HJ}
R. A. Hidalgo and P. Johnson.
Field of Moduli of Generalized Fermat Curves with an application to non-hyperelliptic dessins d'enfants.
{\it Journal of Symbolic Computation} {\bf 71} (2015), 60--72.

\bibitem{HR}
R. A. Hidalgo and S. Reyes-Carocca.
A constructive proof of Weil's Galois descent theorem. 
https://arxiv.org/abs/1203.6294



\bibitem{Huggins-tesis}
B. Huggins.
Fields of Moduli and Fields of Definition of Curves. Ph.D. Thesis, UCLA, 2005.

	
	
\bibitem{Huggins}
B. Huggins. 
Fields of moduli of hyperelliptic curves. 
{\it Math. Res. Lett.} {\bf 14} No.2 (2007), 249--262.


\bibitem{Hurwitz}
A. Hurwitz.
 \"Uber algebraische Gebilde mit eindeutigen Transformationen in sich.
 {\it Math. Ann.} {\bf 41} (1893) 403--442.


\bibitem{Koizumi}
S. Koizumi.
The fields of moduli for polarized abelian varieties and for curves.
{\it Nagoya Math. J.} {\bf 48} (1972) 37--55.


\bibitem{Kontogeorgis}
A. Kontogeorgis. 
Field of moduli versus field of definition for cyclic covers of the projective line.
{\it J. de Theorie des Nombres de Bordeaux} {\bf 21} (2009) 679--692.



\bibitem{LR}
R. Lercier and C. Ritzenthaler.
Hyperelliptic curves and their invariants: geometric, arithmetic and algorithmic aspects. 
{\it J. Algebra} {\bf 372} (2012), 595--636. 

\bibitem{LRRS}
R. Lercier, C. Ritzenthaler, F. Rovetta and J. Sisling.
Parametrizing the moduli space of curves and applications tosmooth plane quartics over finite fields.
{\it LMS J. Comput. Math.} {\bf 17} (Special issue A) (2014) 128--147.

\bibitem{LRS}
R. Lercier, C. Ritzenthaler and J. Sijsling.
Explicit Galois obstruction and descent for hyperelliptic curves with tamely cyclic reduced automorphism group.
{\it Math. Comp.} {\bf 85} (2016), 2011--2045. 


\bibitem{Mestre}
J-F. Mestre.  
Construction de courbes de genre $2$ \`a partir de leurs modules. (French) [Constructing genus-$2$ curves from their moduli] Effective methods in algebraic geometry (Castiglioncello, 1990), 313--334, {\it Progr. Math.} {\bf 94}, Birkh\"auser Boston, Boston, MA, 1991.

\bibitem{KFT}
R. E. Rodr\'{\i}guez and V. Gonz\'alez-Aguilera. 
Fermat's Quartic Curve, Klein's Curve and the Tetrahedron. 
{\it In Contemporary Mathematics} {\bf 201} 1997, 43--62.

\bibitem{Schwarz}
H. A. Schwartz.
\"Uber diejenigen algebraischen Gleichungen zwischen zwei ver\"anderlichen Gr\"o{\ss}en, welche eine schaar rationaler, eindeutig umkehrbarer 
Transformationen in sich selbst zulassen.
{\it Journal f\"ur die reine und angewandte Mathematik} {\bf 87} (1890), 139--145. 


\bibitem{Shimura}
G. Shimura.
On the field of rationality for an abelian variety. 
{\it Nagoya Math. J.} {\bf 45} (1972) 167--178. 



\bibitem{Sil}
J. H. Silverman.
{\it The Arithmetic of Elliptic Curves}.
Graduate Texts in Mathematics {\bf 106}.
Springer, Dordrecht, 2009. xx+513 pp. ISBN: 978-0-387-09493-9




\bibitem{Singerman} 
D. Singerman. 
Finitely Maximal Fuchsian Groups. 
{\it J. London  Math. Soc.} {\bf 6} No.2 (1972) 29--38.

\bibitem{Swinarski}
D. Swinarski.
Equations of Riemann surfaces of genus $4$, $5$ and $6$ with large automorphisms groups.
https://faculty.fordham.edu/dswinarski/publications/SwinarskiEquations.pdf

\bibitem{Weil}
A. Weil.
The field of definition of a variety.
{\it  Amer. J. Math.} {\bf 78} (1956) 509-524.

 \bibitem{Wolfart}
 J. Wolfart.
 $ABC$ for polynomials, dessins d'enfants and uniformization---a survey. Elementare und analytische Zahlentheorie, 313--345, Schr. Wiss. Ges. Johann Wolfgang Goethe Univ. Frankfurt am Main, 20, Franz Steiner Verlag Stuttgart, Stuttgart, 2006.


\end{thebibliography}
\end{document}